\documentclass[12pt,notitlepage]{amsart} 
\usepackage{amscd}

\advance\voffset by -0.4 truein
\advance\hoffset by -0.8 truein
\renewcommand{\:}{\colon}
\newcommand{\QQ}{{\mathbb Q}}
\newcommand{\ol}{\overline}
\newcommand{\wt}{\widetilde}
\newcommand{\cal}{\mathcal}
 \font\tenbi=cmmi14
 \font\sevenbi=cmmi10 \font\fivebi=cmmi7
 \newfam\bifam \def\bi{\fam\bifam} \textfont\bifam=\tenbi
 \scriptfont\bifam=\sevenbi \scriptscriptfont\bifam=\fivebi
 \mathchardef\variablemega="7121 \def\w{{\bi\variablemega}}

\newtheorem{theorem}{Theorem}[section]
\newtheorem{proposition}[theorem]{Proposition}
\theoremstyle{definition}  
\newtheorem{claim}[theorem]{}

\font\smallrm=cmr8
\font\midrm=cmr10
\font\midbf=cmbx10

\begin{document}

\

\vskip2cm

\author[{\smallrm Eduardo Esteves -- Letterio Gatto}]
{\vskip0.5cm Eduardo Esteves -- Letterio Gatto}
\thanks{First author supported by PRONEX, Conv\^enio 41/96/0883/00, CNPq, 
Proc. 300004/95-8, and FAPERJ, Proc. E-26/170.418/2000-APQ1.}
\thanks{Second author supported by CNPq, Proc. 300680/99-6, GNSAGA-CNR 
and MURST, Progetto Nazionale ``Geometria Algebrica, Algebra Commutativa e
Aspetti Computazionali'', coordinatore Claudio Pedrini.}
\title[{\smallrm On a relation by Cornalba and Harris}]
{A geometric interpretation and a new proof 
of a relation by Cornalba and Harris}
%
%
\maketitle

\vskip0.5cm

\begin{center}
\begin{minipage}[r]{5.5in}{
\baselineskip=12pt
\noindent
{\midbf Abstract.} \midrm 
In the 80's M.~Cornalba and J.~Harris discovered a 
relation among the Hodge class and the boundary classes in the 
Picard group with rational coefficients of the moduli space 
of stable, hyperelliptic 
curves. They proved the relation by computing degrees of 
the classes involved for suitable one-parameter families. In the 
present article we show that their relation can be obtained as the class of 
an appropriate, geometrically meaningful empty set, thus conforming with 
C.~Faber's general philosophy 
to finding relations among tautological classes in the Chow ring of 
the moduli space of curves. 
The empty set we consider is the closure of the locus of smooth, 
hyperelliptic curves having a special ramification point.

\vskip0.2cm

}
\end{minipage}
\end{center}

\vskip1.5cm

\section{Introduction}
\label{Sect1}

\begin{claim}\emph{Mumford's relation and Faber's idea.}
Let $\ol M_g$ be the moduli space of 
(Deligne-Mumford) stable curves of genus $g$, and 
$\text{\rm Pic}(\ol M_g)$ its Picard group; see \cite{DM} and \cite{HaMo}.
Natural classes of $\text{\rm Pic}(\ol M_g)\otimes\QQ$ are the 
Hodge class $\lambda$ and the boundary classes 
$\delta_0,\dots,\delta_{[g/2]}$; see \cite{AC}. 
It is a fundamental result by J. Harer, 
E. Arbarello and M. Cornalba 
that the above classes freely generate 
$\text{\rm Pic}(\ol M_g)\otimes\QQ$ for each $g\geq 3$; see \cite{Harer} 
and \cite{AC}. However, for $g=2$ there is the relation below, 
proved by D. Mumford in his fundamental paper~\cite{MumfEnum},
\begin{equation}
10\lambda-\delta_0-2\delta_1=0.\label{eq:Mumrel}
\end{equation}

Mumford's relation~(\ref{eq:Mumrel}) was recovered by S. Diaz and 
F. Cukierman in \cite{Diazexc} and \cite{cukie}, though 
they do not observe it there. In fact, Diaz considers 
the locus $D_{g-1,g-1}\subseteq\ol M_g$ 
of smooth curves having a (Weierstrass) point whose first non-zero 
non-gap is 
at most $g-1$, and computes the class of its closure in 
$\text{\rm Pic}(\ol M_g)\otimes\QQ$ (\cite{Diazexc}, Thm.~(7.4), p. 40). 
Similarly, Cukierman computes
the class of the closure of the locus $\cal E\subseteq\ol M_g$ 
of smooth curves having a global holomorphic 1-form vanishing 
at a (Weierstrass) point with order at least 
$g+1$ (\cite{cukie}, Eq.~(5.5), p.~344 and Rmk.~(b) thereafter).
Both Diaz and 
Cukierman claim their formulas hold for $g\geq 3$. However, we observed 
that their formulas hold for $g=2$ as well, and give Mumford's 
relation~(\ref{eq:Mumrel}). Indeed, for $g=2$ both $D_{1,1}$ and $\cal E$ 
are empty, and 
so the formulas given by Diaz and Cukierman 
express the class of 
the empty set, 0, as a linear combination of $\lambda$, $\delta_0$ and 
$\delta_1$. In both cases, this linear combination is a non-zero multiple 
of that in~(\ref{eq:Mumrel}).

Therefore, Mumford's relation~(\ref{eq:Mumrel}) fits perfectly in 
C.~Faber's 
philosophy to finding relations among tautological classes in Chow 
rings of moduli spaces of curves by writing down classes of 
``appropriate'' empty 
sets; see \cite{Faber3}. The beauty in Faber's idea is that relations 
discovered by his approach are endowed with geometric interpretations. 

Now, let $H_g\subseteq\ol M_g$ be the locus of smooth, hyperelliptic 
curves, and $\ol H_g\subseteq\ol M_g$ its closure. A stable curve 
represented in $\ol H_g$ is also called hyperelliptic. 
Of course, each smooth curve of 
genus 2 is hyperelliptic, hence $\ol H_2=\ol M_2$. So we may view 
Mumford's relation as one in $\text{Pic}(\ol H_2)\otimes\QQ$. 
In general, M. Cornalba and J. Harris proved in \cite{CH} 
the theorem below; see also Subsection~6C of \cite{HaMo}.

\begin{theorem} {\rm (Cornalba-Harris)} The following relation 
holds in $\text{\rm Pic}(\ol H_g)\otimes\QQ$:
\begin{equation}
(8g+4)\lambda-g\xi_0-\sum_{i=1}^{[(g-1)/2]}2(i+1)(g-i)\xi_i-
\sum_{j=1}^{[g/2]}4j(g-j)\delta_j=0.\label{eq:CH}
\end{equation}
\end{theorem}

In the above statement, $\xi_0,\dots,\xi_{[(g-1)/2]}$ are boundary classes 
satisfying the relation below; see Sections~\ref{Sect3}~and~\ref{Sect4}. 
\begin{equation}
\delta_0-\xi_0-2(\xi_1+\cdots+\xi_{[(g-1)/2]})=0.\label{eq:deltaxi}
\end{equation}
For $g=2$, Cornalba's and Harris' relation~(\ref{eq:CH}) coincides with 
Mumford's~(\ref{eq:Mumrel}).

Cornalba and Harris found~(\ref{eq:CH}) by computing the degree of 
$\lambda$ for one-parameter families of admissible covers. Now, our 
observation above suggests that, like Mumford's, also Cornalba's and Harris' 
relation might be obtained by expressing the 
class of an ``appropriate'' empty set, thus shedding light on its geometric 
significance. As in the works by Diaz and Cukierman, the ``appropriate'' 
empty set might be the locus of curves having a certain type of 
Weierstrass point. Indeed, the goal of the present article is to show 
that \emph{a non-zero multiple of 
the left-hand side of~{\rm (\ref{eq:CH})} 
expresses the class in $\text{\rm Pic}(\ol H_g)\otimes\QQ$ of 
the closure $\ol W_{g,2}$ of the locus $W_{g,2}\subseteq\ol H_g$ of smooth, 
hyperelliptic curves 
for which the unique $g^1_2$ has a special ramification point}.
\end{claim}

\begin{claim}\emph{An overview.} 
As Cornalba and Harris in \cite{CH}, we will rather work with the moduli 
stack $\ol{\cal H}_g$ of stable, hyperelliptic curves of genus $g$, and 
show relation~{\rm (\ref{eq:CH})} in 
$\text{\rm Pic}(\ol{\cal H}_g)\otimes\QQ$. In other words, we will 
consider a family $\pi\:C\to S$ of stable, hyperelliptic curves and show 
that the induced classes $\lambda_\pi$, 
$\xi_{1,\pi},\dots,\xi_{[(g-1)/2],\pi}$ and 
$\delta_{1,\pi},\dots,\delta_{[g/2],\pi}$ in 
$\text{\rm Pic}(S)\otimes\QQ$ satisfy the relation
\begin{equation}
(8g+4)\lambda_\pi-g\xi_{0,\pi}-\sum_{i=1}^{[(g-1)/2]}2(i+1)(g-i)\xi_{i,\pi}-
\sum_{j=1}^{[g/2]}4j(g-j)\delta_{j,\pi}=0.\label{eq:CHpi}
\end{equation}
We will see that it is sufficient to consider 
a particular family $\pi\:C\to S$ where $S$ is a smooth, projective curve, 
the general fiber of $\pi$ is smooth and admits a $g^1_2$, and the singular 
fibers are ``general enough''; see Section~\ref{Sect4}. 
So let us assume $\pi$ is this family.

Let $\wt\pi\:\wt C\to S$ be the semi-stable 
reduction of $\pi$. In Section~\ref{Sect5} we consider the 
wronskian determinant $w$ of a certain relative $g^1_2$ on $\wt\pi$, and use 
the theory of limit linear series developed in~\cite{EST2} to compute
the orders of vanishing of $w$ on the irreducible components of the 
singular fibers of $\wt\pi$. Subtracting these irreducible components with 
the computed multiplicities, we obtain a section $\ol w$ of a line 
bundle on $\wt C$ cutting out a relative Cartier divisor $\ol W$ over $S$. 
Now, as in~\cite{GP} (see also~\cite{Gatto2}), 
we consider in Section~\ref{Sect6} the degeneracy locus $Z$ of the 
derivative $D{\ol w}$ of $\ol w$. This locus is empty, as we show 
that $\ol W$ is unramified over $S$ for the particular family $\pi$ we 
consider. Since $Z$ has the 
expected dimension, the class of $[Z]$ coincides with the 
degeneracy class of $D{\ol w}$, what allows us to compute an expression 
for $\wt\pi_*[Z]$ in terms of $\lambda_\pi$, 
$\xi_{1,\pi},\dots,\xi_{[(g-1)/2],\pi}$ and 
$\delta_{1,\pi},\dots,\delta_{[g/2],\pi}$.
This expression 
is the class in $\text{Pic}(S)\otimes\QQ$ of the closure of the set of 
points $s\in S$ for which the 
fiber $C_s$ is smooth and has a special ramification point for its 
unique $g^1_2$. This set is empty, and the expression is 
a multiple of the left-hand side of~(\ref{eq:CHpi}). Thus~(\ref{eq:CHpi}) 
is proved.

In Section~\ref{Sect2} we exemplify our method by 
showing~(\ref{eq:Mumrel}). This is simpler than showing~(\ref{eq:CH}), 
because all the singular curves we have to deal with are either 
irreducible or of compact type, and hence we may apply results by Diaz and 
Cukierman in \cite{Diazexc} and \cite{cukie}. So, not only are the 
technical tools kept to a minimum, but also are available in the literature. 
To conform to the technical tools available, the proof in 
Section~\ref{Sect2} is just slightly different from the one laid out above.
\end{claim}

\begin{claim}\emph{Acknowledgements.} This collaboration began during 
the second author's visit to IMPA in February, 2000, and was 
concluded during the first author's visit to UFPe in June of the same 
year. Both institutions are warmly thanked.
\end{claim}

\section{Mumford's relation}
\label{Sect2}

\begin{claim}\emph{Reduction.} 
Let $\pi\:C\to S$ be a family of stable curves of genus $g$. Denote by 
$\w_\pi$ its relative dualizing sheaf. The direct image 
$E_\pi:=\pi_*\w_\pi$ is called the \emph{Hodge bundle} of $\pi$. Let 
$\lambda_\pi:=c_1(E_\pi)$ in $\text{Pic}(S)$. 
For each invertible 
sheaf ${\cal L}$ on $C$, let $J^1_\pi(\cal L)$ be the 
\emph{sheaf of relative jets}, or principal parts of order 1 of $\cal L$; 
see Section~2 of \cite{Gatto2} or \cite{EST1}. As for families of smooth 
curves, there is a natural exact sequence,
\begin{equation}
0\to\cal L\otimes\w_\pi\to J^1_\pi(\cal L)\to\cal L\to 0.\label{eq:short}
\end{equation}

If $S$ and the general fiber of $\pi$ are smooth curves 
then the singularities of the 
total space $C$ of $\pi$ are concentrated at the nodes of the singular 
fibers and are mild. More precisely, 
the surface $C$ has at each node $P$ of each 
singular fiber of $\pi$ a singularity of type $A_k$ for a certain 
integer $k$, called the 
\emph{singularity type of $P$ in $C$}. Then $C$
is smooth if and only if the singularity type in $C$ of every node of 
every singular fiber of $\pi$ is 1. 

Using level structures, one can show that there is a family of stable 
curves over a smooth base $\Omega$ represented by a finite and surjective 
map $\Omega\to\ol M_g$; see 
\cite{L} or \cite{B}. Hence, to show a relation among classes in 
$\text{\rm Pic}(\ol M_g)\otimes\QQ$, it is enough to prove that for 
a ``sufficiently general'' one-parameter family of stable curves 
the induced classes satisfy the corresponding relation. 

Assume $g=2$ for the rest of Section~\ref{Sect2}. Let 
$\pi\: C\to S$ be an one-parameter family of stable curves of genus 2. 
Assume $\pi$ is ``sufficiently general''. More precisely, assume $S$ 
and the general fiber of $\pi$ are smooth curves, and each of the singular 
fibers of $\pi$ is a general uninodal curve; in particular, each fiber of 
$\pi$ has at most two irreducible components. To 
prove~(\ref{eq:Mumrel}) we need only show that the classes 
$\lambda_\pi,\delta_{0,\pi},\delta_{1,\pi}$ in $\text{Pic}(S)\otimes\QQ$ 
induced by $\lambda,\delta_0,\delta_1$ satisfy
\begin{equation}
10\lambda_\pi-\delta_{0,\pi}-2\delta_{1,\pi}=0.
\label{eq:Mumrelpi}
\end{equation}

For each $s\in S$ such that 
$C_s$ is singular, let $k_s$ be the singularity type in $C$ of the unique 
node of $C_s$. 
For each $i=0,1$, let $\Delta_{i,\pi}$ be the set of points $s\in S$ for 
which the fiber $C_s$ is singular with $i+1$ irreducible components. 
Then $\delta_{i,\pi}$ is the divisor class of points 
$s\in\Delta_{i,\pi}$ counted each with multiplicity $k_s$; 
see~Section~3D~of~\cite{HaMo}. 
\end{claim}

\begin{claim}\emph{Computation.} 
Let $w$ be the section of the line bundle 
$\w_\pi^3\otimes\pi^*\wedge^2E_\pi^*$ 
obtained by taking determinants in 
the natural map of vector bundles $\pi^*E_\pi\to J^1_\pi(\w_\pi)$. 
The zero scheme $W$ of $w$ 
cuts each smooth fiber of $\pi$ in the scheme of Weierstrass points 
of that fiber. Using the local computations by Cukierman for 
Prop.~(2.0.8) on p.~325 of~\cite{cukie}, 
we see that $w$ vanishes 
along each irreducible component of $C_s$ for each $s\in\Delta_{1,\pi}$ with 
order $k_s$. Let $\ol w$ be the induced section of 
$\w_\pi^3\otimes\pi^*(\wedge^2E_\pi^*(-D))$, where 
$D:=\sum_{s\in\Delta_1}k_s[s]$. Let $\ol W$ be the zero scheme of 
$\ol w$. Then $\ol W$ is a 
relative Cartier divisor on $C$ over $S$. By considerations of degree, 
away from possible nodes the divisor $\ol W$ cuts each 
fiber $C_s$ for each $s\in S$ transversally. 
By Theorem~(A2.4) on page~64 of \cite{Diazexc}, for each 
$s\in\Delta_{1,\pi}$ the divisor $\ol W$ cuts the fiber $C_s$ away from 
the node. In contrast, by Theorem~(A2.1) on page~60 of \cite{Diazexc}, 
the divisor $\ol W$ cuts the fiber $C_s$ for each $s\in\Delta_{0,\pi}$ 
with multiplicity 2 at the node.  

Let $Z$ be the degeneracy locus of the derivative $D\ol w$ of $\ol w$. 
Here $D\ol w$ is the section of 
\begin{equation}
J^1_\pi(\w_\pi^3\otimes\pi^*(\wedge^2E_\pi^*(-D)))\label{eq:J1}
\end{equation}
induced by $\ol w$ (see \cite{Gatto2} or \cite{GP} for this 
construction). Since $\ol W$ is a Cartier divisor cutting the fiber $C_s$ 
for each $s\in\Delta_{0,\pi}$ with multiplicity 2 at the node, a local 
computation shows that the length of $Z$ at the node is $k_s$. As 
$Z$ is supported in the union of these nodes, we have 
$\pi_*[Z]=\delta_{0,\pi}$. 
Now, as the degeneracy locus of $D\ol w$ has the expected dimension zero, 
the class of $[Z]$ is the second Chern class of the bundle in~(\ref{eq:J1}). 
So
\begin{align*}
\pi_*[Z]=&
\pi_*[(3c_1(\w_\pi)-\pi^*(\lambda_\pi+\delta_{\pi,1}))
(4c_1(\w_\pi)-\pi^*(\lambda_\pi+\delta_{\pi,1}))]\\
=&\pi_*[12c_1(\w_\pi)^2-7c_1(\w_\pi)\pi^*(\lambda_\pi+\delta_{1,\pi})].
\end{align*}
By Grothendieck-Riemann-Roch, 
$\pi_*(c_1(\w_\pi)^2)=12\lambda_\pi-\delta_{0,\pi}-\delta_{1,\pi}$; 
see Section 3E~of~\cite{HaMo} or 
Theorem~5.10~on~page~100~of~\cite{MumfStab}. By 
the projection formula, $\pi_*(c_1(\w_\pi)\pi^*\mu)=2\mu$ for 
every $\mu\in\text{Pic}(S)$. Hence,
\[
\pi_*[Z]=12(12\lambda_\pi-\delta_{0,\pi}-\delta_{1,\pi})
-14(\lambda_\pi+\delta_{1,\pi})=130\lambda_\pi-12\delta_{0,\pi}-
26\delta_{\pi,1},
\]
and thus the difference $\pi_*[Z]-\delta_{0,\pi}$ is a multiple of the 
left-hand side of~(\ref{eq:Mumrelpi}). Since $\pi_*[Z]=\delta_{0,\pi}$ 
we get~(\ref{eq:Mumrelpi}), thus proving 
Mumford's relation~(\ref{eq:Mumrel}).
\end{claim}

\section{The Picard group}
\label{Sect3}

\begin{claim}
\emph{The Picard group of $\ol M_g$.} (See \cite{AC} or \cite{HaMo}.)
Let $\Delta_0\subseteq\ol M_g$ be the closure of the locus of irreducible, 
uninodal, stable curves. For each $i=1,\dots,[g/2]$, 
let $\Delta_i\subseteq\ol M_g$ be the closure of the locus of uninodal, 
stable curves with two irreducible components, one of them of genus $i$. 
Then $\Delta_0,\dots,\Delta_{[g/2]}$ are divisors of $\ol M_g$ covering 
$\ol M_g-M_g$. 

Let $\ol{\cal M}_g$ be the moduli stack of stable curves and 
$\text{Pic}(\ol{\cal M}_g)$ its group of divisor classes. Each divisor 
class on $\ol M_g$ defines, by pullback, one on $\ol{\cal M}_g$. 
So, there is a well-defined homomorphism 
$\Psi\:\text{Pic}(\ol M_g)\to\text{Pic}(\ol{\cal M}_g)$. Moreover, one can 
show that $\Psi\otimes 1_{\QQ}$ is an isomorphism.

There is a class $\lambda\in\text{Pic}(\ol{\cal M}_g)$ whose pullback 
to $\text{Pic}(S)$ for each family $\pi\:C\to S$ of stable curves is 
the class $\lambda_\pi$ defined in Subsection~2.1. We call $\lambda$ 
the \emph{Hodge class}. In addition, as the universal 
deformation space of a stable curve is smooth, 
for each $i=0,\dots,[g/2]$ the divisor $\Delta_i$ induces a 
divisor class $\delta_i$ in $\text{Pic}(\ol{\cal M}_g)$. 
We call $\delta_0,\dots,\delta_{[g/2]}$ the 
\emph{boundary classes}. Denote also by 
$\lambda,\delta_0,\dots\delta_{[g/2]}$ the corresponding 
classes in $\text{Pic}(\ol M_g)\otimes\QQ$ under the isomorphism 
$\Psi\otimes 1_{\QQ}$. 
\end{claim}

\begin{claim}
\emph{The Picard group of $\ol H_g$.} (See \cite{CH} or Section~6C of 
\cite{HaMo}.) The locus $\ol H_g$ is the quotient of the Hurwitz scheme 
of admissible covers of degree 2 by the action of the 
symmetric group on $2g+2$ letters. Moreover, the latter scheme is the 
same thing as the moduli space $\ol M_{0,2g+2}$ of stable $(2g+2)$-pointed 
rational curves. The moduli space $\ol M_{0,2g+2}$ is fine, smooth and 
projective.

The above description of $\ol H_g$ implies that $\ol H_g-H_g$ has exactly 
$g$ irreducible components $\Xi_0,\dots,\Xi_{[(g-1)/2]}$ 
and $\Theta_1,\dots,\Theta_{[g/2]}$, all of them divisors of 
$\ol H_g$, given as follows:
\begin{enumerate}
\item The component $\Xi_0$ is the closure of the locus of curves 
obtained from a smooth, hyperelliptic curve of genus $g-1$ by 
identifying two distinct points conjugated by the hyperelliptic involution.
\item For each $i=0,\dots,[(g-1)/2]$, the component $\Xi_i$ is the 
closure of the locus of curves with just two irreducible 
components of genera $i$ and $g-1-i$ which are smooth, hyperelliptic, 
and intersect transversally at just two distinct points conjugated 
by the hyperelliptic involution on each curve.
\item For each $j=1,\dots,[g/2]$, the component $\Theta_j$ 
is the closure of the locus of curves with just two irreducible 
components of genera $j$ and $g-j$ which are smooth, hyperelliptic, 
and intersect transversally at a unique point fixed 
by the hyperelliptic involution on each curve. 
\end{enumerate}

As for stable curves, there is a natural isomorphism 
$\Phi\:\text{Pic}(\ol H_g)\to\text{Pic}(\ol{\cal H}_g)$, where 
$\ol{\cal H}_g$ is the moduli stack of stable, hyperelliptic curves of 
genus $g$. As pointed out on page~469 of \cite{CH}, since the 
universal deformation space of a stable, 
hyperelliptic curve within hyperelliptic curves is smooth, 
for each $i=0,\dots,[(g-1)/2]$ the divisor $\Xi_i$ induces a 
divisor class $\xi_i$ in $\text{\rm Pic}(\ol{\cal H}_g)$. 
Denote also by $\xi_i$ the corresponding 
rational divisor class in $\ol H_g$ under the isomorphism 
$\Phi$. Though improperly, we shall also denote 
by $\lambda,\delta_0,\dots,\delta_{[g/2]}$ the pullbacks to $\ol H_g$ or 
$\ol{\cal H}_g$ of the corresponding classes in $\ol M_g$ or 
$\ol{\cal M}_g$. Then relation~(\ref{eq:deltaxi}) is simply a restatement 
of the identity~(4.6) on page 469 of \cite{CH}. Relation~(\ref{eq:deltaxi}) 
is also recovered in Section~\ref{Sect4} below.
\end{claim}

\section{Reduction}
\label{Sect4}

The natural map $\ol M_{0,2g+2}\to\ol H_g$ is finite and surjective. So, 
the pullback homomorphism,
\[
\text{Pic}(\ol H_g)\otimes\QQ\to\text{Pic}(\ol M_{0,2g+2})\otimes\QQ,
\]
is injective.
Since $\ol M_{0,2g+2}$ is smooth and projective, there is a smooth 
curve $T\subseteq\ol M_{2g+2}$ such that the left-hand sides 
of~(\ref{eq:CH}) and~(\ref{eq:deltaxi}) pull back to zero on 
$\ol M_{0,2g+2}$ only if they pull back to zero on $T$. 
Since $\ol M_{0,2g+2}$ is a fine moduli 
space, there are a family $R/T$ of rational, nodal curves and 
sections $\zeta_1,\dots,\zeta_{2g+2}$ of the smooth locus of 
$R/T$ representing the 
inclusion $T\to\ol M_{0,2g+2}$. The curve $T$ may be chosen 
``general enough'' so that $R$ is smooth and each singular fiber of 
$R/T$ is uninodal.

Let $T_0$ (resp. $T_1$, resp. $T_2$) be the subset of points 
$t\in T$ for which $R_t$ is 
singular and one of its irreducible components contains two 
(resp. an odd number, resp. an even number) of the points 
$\zeta_1(t),\dots,\zeta_{2g+2}(t)$. 
We can construct a smooth curve $S$ and a finite 
covering $\gamma\:S\to T$ of degree 2 ramified along $T_1$ but unramified 
along $T_0\cup T_2$; see Lecture~3 in \cite{EV} for the theory of 
cyclic coverings. 
Put $Q:=R\times_T S$ and 
$S_i:=\gamma^{-1}(T_i)$ for each $i=0,1,2$. Let $B$ be 
the blow-up of $Q$ along the nodes of the fibers $Q_s$ 
for all $s\in S_1$. Then $B$ is smooth. For each $s\in S_1$, let 
$E_s$ be the rational curve on $B$ contracting to a point on 
the fiber $Q_s$. For each $i=1,\dots,2g+2$, as 
$\zeta_i$ is a section of $Y/T$ through its smooth locus, 
$\zeta_i\times 1_S\:S\to Q$ factors through the blow-up map $B\to Q$; 
let $\tau_i$ denote this section of $B/S$. Note that 
$\tau_i(S)\cap E_s=\emptyset$ for each $i=1,\dots,2g+2$ and each 
$s\in S_1$. So the divisor $D$ defined below is reduced.
\[
D:=\sum_{i=0}^{2g+2}\tau_i(S)+\sum_{s\in S_1}E_s. 
\]
Moreover, $D$ is smooth, and its restriction 
to each irreducible component of each fiber of $B/S$ has even degree. 
So, we can construct a smooth surface $M$ and a finite covering 
$M\to B$ of degree 2 ramified along $D$. 
For each $s\in S_1$, let $F_s$ be the rational curve on $M$ 
lying over $E_s$. Let $\wt C$ be the surface 
obtained from $M$ 
by contracting $F_s$ for all $s\in S_1$, and $\wt\pi\:\wt C\to S$ 
the induced map. Since $F_s$ has self-intersection $-1$, the surface $\wt C$ 
is smooth. Finally, 
let $C$ be the surface obtained from $\wt C$ 
by contracting all rational curves contained in the fibers over $S_0$. 
Then the 
induced map $\pi\:C\to S$ is a family of 
stable, hyperelliptic curves representing the composition 
$\beta\:S\to T\to\ol M_{0,2g+2}\to\ol H_g$. 

For each $i=1,\dots,2g+2$, the section $\tau_i$ factors through the 
covering $M\to B$ because the latter is ramified along $\tau_i(S)$; let 
$\sigma_i$ denote this section of $M/S$. Then $\sigma_i$ induces a section 
of $\wt\pi$; let 
$\Sigma_i$ denote the image of this section in $\wt C$. It follows from 
our construction that $\Sigma_i$ is contained in the smooth locus of 
$\wt C/S$. 

Put $\Sigma:=\Sigma_1$ and $\cal L:=\mathcal O_{\wt C}(2\Sigma)$. 
Let $\ol W:=\Sigma_1+\cdots+\Sigma_{2g+2}$. 
Then $\ol W$ is flat over $S$, and intersects each 
smooth fiber of $\wt\pi$ in the ramification scheme of the 
complete linear system of sections of the restriction of $\cal L$ 
to that fiber. 

From now on to the end of the article, fix $\pi$ and $\wt\pi$ as above. 
Since $\gamma\:S\to T$ is finite and surjective, 
the pullback map 
$\gamma^*\:\text{Pic}(T)\otimes\QQ\to\text{Pic}(S)\otimes\QQ$ is 
injective. Hence, (\ref{eq:CH}) holds only if~(\ref{eq:CHpi}) holds, 
whereas~(\ref{eq:deltaxi}) holds only if 
\begin{equation}
\delta_{0,\pi}-\xi_{0,\pi}-
2(\xi_{1,\pi}+\cdots+\xi_{[(g-1)/2],\pi})=0.\label{eq:deltaxipi}
\end{equation}

For each $i=0,\dots,[(g-1)/2]$ and each $j=1,\dots,[g/2]$, denote
$\Xi_{i,\pi}:=\beta^{-1}(\Xi_i)$ and 
$\Delta_{j,\pi}:=\beta^{-1}(\Theta_j)$ with their reduced 
induced structures. Let
\[
\Delta_\pi:=\Delta_{1,\pi}\cup\dots\cup\Delta_{[g/2],\pi}\quad{\rm and}\quad 
\Xi_\pi:=\Xi_{0,\pi}\cup\dots\cup\Xi_{[(g-1)/2],\pi}.
\]
Put $\Gamma_\pi:=\Delta_\pi\cup\Xi_\pi$. Then 
$\Gamma_\pi$ consists of the points $s\in S$ for which 
the fiber $C_s$ is singular. By construction, $C$ is smooth everywhere 
but at the nodes of the fibers of $\pi$ over $\Xi_{0,\pi}$, 
where $C$ has singularities of type $A_2$. It follows that
\begin{align}
\xi_{0,\pi}=&2[\Xi_{0,\pi}]\quad\text{and}\quad\xi_{i,\pi}=[\Xi_{i,\pi}]
\  \text{for $i=1,\dots,[(g-1)/2]$},\label{eq:DeltaXipi1}\\
\delta_{0,\pi}=&2[\Xi_\pi]\quad\text{and}\quad
\delta_{j,\pi}=[\Delta_{j,\pi}]\  \text{for $j=1,\dots,[g/2]$}.
\label{eq:DeltaXipi2}
\end{align}
Relation~(\ref{eq:deltaxipi}) follows immediately from 
the above expressions. Hence~(\ref{eq:deltaxi}) holds. 

\section{Limit linear systems and ramification divisors}
\label{Sect5}

\begin{claim}\emph{Limit linear systems.} For each 
$s\in\Gamma_\pi$, let $X_s$ and $Y_s$ be the irreducible components of 
$\wt C_s$, chosen such that $\Sigma_s\subseteq X_s$. Let 
$a_s$ and $b_s$ be the genera of $X_s$ and $Y_s$, respectively. 

By \cite{EST1} or \cite{Ran}, there are unique (mod ${\rm Pic}(S)$) 
invertible sheaves $\cal I$ and $\cal K$ on $\wt C$ which agree with 
$\cal L$ on the generic fiber of $\wt\pi$ and satisfy the 
properties below for each $s\in\Gamma_\pi$:
\begin{enumerate}
\item The natural maps $\wt\pi_*\cal I(s)\to H^0(\cal I|X_s)$ and 
$\wt\pi_*\cal K(s)\to H^0(\cal K|Y_s)$ are injective.
\item The natural maps $\wt\pi_*\cal I(s)\to H^0(\cal I|Y_s)$ 
and $\wt\pi_*\cal K(s)\to H^0(\cal K|X_s)$ are non-zero.
\end{enumerate}
The sheaves $\wt\pi_*\cal I$ and $\wt\pi_*\cal K$ are generically of 
rank 2. Moreover, they are torsion-free because $\wt\pi$ is flat. So, 
they are locally free of rank 2 because $S$ is smooth. 

Let $s\in\Gamma_\pi$. 
By Properties~1~and~2 above, the restriction of 
$\cal I$ to either $X_s$ or $Y_s$ has non-negative degree. Suppose 
first that $X_s$ is not rational. Then $\deg\cal I|X_s\geq 2$ because 
$\dim\wt\pi_*\cal I(s)=2$. Now, $\cal I$ has total degree 2, as 
it agrees with $\cal L$ on the generic fiber of $\wt\pi$. It follows that 
$\deg\cal I|X_s=2$ and $\deg\cal I|Y_s=0$.

Now, suppose $X_s$ is rational. Then $Y_s$ is not rational, 
and hence $\deg\cal K|X_s=0$ by analogy with the case above. 
Since $\cal I$ and $\cal K$ 
agree on the generic fiber of $\wt\pi$, there is an integer 
$\ell$ such that $\cal I$ agrees with $\cal K(\ell Y_s)$ on a neighborhood 
of $\wt C_s$. Since $\deg\cal K|X_s=0$, it follows from Property~1 
that $\ell>0$. Now, since $X_s$ is rational, 
we have $s\in\Xi_{0,\pi}$ and, in particular, $X_s\cdot Y_s=2$. So 
$\deg\cal I|X_s\geq 2$ because $\ell>0$. As before, 
$\deg\cal I|X_s=2$ and $\deg\cal I|Y_s=0$. 

By analogy, $\deg\cal K|X_s=0$ and $\deg\cal K|Y_s=2$.

Since $\cal L$ and $\cal I$ agree on the generic fiber of $\wt\pi$, and 
restrict to sheaves of equal degree on 
each irreducible component of each singular fiber, $\cal L$ and $\cal I$ 
agree mod $\text{Pic}(S)$. We may assume that $\cal I=\cal L$. 
By the same reason, we may assume that 
\begin{equation}
\cal K=\cal L(-\sum_{s\in\Xi_\pi}Y_s-2\sum_{s\in\Delta_\pi}Y_s).\label{eq:N}
\end{equation}

If $s\in S$ is such that $C_s$ is smooth, 
then $\cal L_s^{g-1}=\w_{\wt\pi,s}$. Now, for each $s\in\Gamma_\pi$,
\[
\deg\w_{\wt\pi}|Y_s=
\begin{cases}
2b_s&\text{if $s\in\Xi_\pi$},\\
2b_s-1&\text{if $s\in\Delta_\pi$}.
\end{cases}
\]
Since $\deg\cal L|Y_s=0$, by considerations of degree, 
\begin{equation}
{\cal L}^{g-1}=\w_{\wt\pi}(\sum_{s\in\Xi_\pi}b_sY_s+
\sum_{s\in\Delta_\pi}(2b_s-1)Y_s)
\quad\text{(mod $\text{Pic}(S)$)}.\label{eq:Lg-1}
\end{equation}
\end{claim}

\begin{claim}\emph{Ramification divisors.} Let $W$ be the 
ramification divisor associated to 
$\cal I$ and $W'$ that associated to 
$\cal K$. In other words, $W$ and $W'$ are the degeneracy 
schemes of the natural maps $\wt\pi^*\wt\pi_*\cal I\to J^1_{\wt\pi}(\cal I)$ 
and $\wt\pi^*\wt\pi_*\cal K\to J^1_{\wt\pi}(\cal K)$, respectively. 
Then $\ol W$ agrees with $W$ and $W'$ away from the singular fibers of 
$\wt\pi$. Moreover, by Property~1 of Subsection~5.1, 
neither $X_s$ is contained 
in the support of $W$, nor $Y_s$ is contained in the support of $W'$. 
\end{claim}

\begin{proposition} The following expression holds:
\begin{equation}
\ol W=W-\sum_{s\in\Xi_\pi}Y_s-2\sum_{s\in\Delta_\pi}Y_s.\label{eq:olW}
\end{equation}
\end{proposition}

\begin{proof} Let
\begin{equation}
\cal J:=\cal I(-\sum_{s\in\Gamma_\pi}Y_s).\label{eq:P}
\end{equation}
Let $w$, $w^\dagger$ and $w'$ be the respective 
determinants of the natural maps,
\[
{\wt\pi}^*\wt\pi_*\cal I\longrightarrow J^1_{\wt\pi}(\cal I),\quad 
{\wt\pi}^*\wt\pi_*\cal J\longrightarrow J^1_{\wt\pi}(\cal J)\quad\text{and}
\quad{\wt\pi}^*\wt\pi_*\cal K\longrightarrow J^1_{\wt\pi}(\cal K).
\]
The zero scheme of $w$ is $W$ and that of $w'$ is $W'$. 
Let $W^\dagger$ be the zero scheme of $w^\dagger$. 

As observed in Subsection~5.1, we have $\cal I=\cal L$. 
By~(\ref{eq:N})~and~(\ref{eq:P}), there is a sequence of inclusions 
$\cal K\hookrightarrow\cal J\hookrightarrow\cal I$, which induces the 
commutative diagram below.
\[
\begin{CD}
{\wt\pi}^*\det\wt\pi_*\cal K @>>> {\wt\pi}^*\det\wt\pi_*\cal J @>>> 
{\wt\pi}^*\det\wt\pi_*\cal I\\
@Vw'VV @Vw^\dagger VV @VwVV\\
\w_{\wt\pi}\otimes{\cal K}^2 @>2\sum_{s\in\Delta_\pi}Y_s>> 
\w_{\wt\pi}\otimes{\cal J}^2 @>2\sum_{s\in\Gamma_\pi}Y_s>> 
\w_{\wt\pi}\otimes{\cal I}^2.
\end{CD}
\]
The horizontal maps above are injective, so we will view them as 
inclusions. By Property~2 of 
Subsection~5.1, there is a local section of $\wt\pi_*\cal I$ that is not 
zero on $Y_s$. Now, $\wt\pi_*\cal I$ has rank 2. 
Since $\deg\cal I|Y_s=0$, there is also 
a local section of $\wt\pi_*\cal I$ that 
vanishes on $Y_s$. So
\begin{equation}
\det\wt\pi_*\cal J=(\det\wt\pi_*\cal I)(-\Gamma_\pi).\label{eq:w}
\end{equation}
Now, by~(\ref{eq:N})~and~(\ref{eq:P}),
\[
\cal J\otimes{\wt\pi}^*{\cal O}_S(-\Delta_\pi)=
\cal K(-\sum_{s\in\Delta_\pi}X_s).
\]
Hence, by analogy with~(\ref{eq:w}),
\begin{equation}
\det\wt\pi_*\cal K=(\det\wt\pi_*\cal J)(-\Delta_\pi).\label{eq:w'}
\end{equation}

Since the two squares in the above diagram are commutative, 
(\ref{eq:w}) and (\ref{eq:w'}) imply
\[
W+\sum_{s\in\Gamma_\pi}X_s=W^\dagger+\sum_{s\in\Gamma_\pi}Y_s
\quad\text{and}\quad 
W^\dagger+\sum_{s\in\Delta_\pi}X_s=W'+\sum_{s\in\Delta_\pi}Y_s,
\]
and hence 
\begin{equation}
W-\sum_{s\in\Gamma_\pi}Y_s-\sum_{s\in\Delta_\pi}Y_s=
W'-\sum_{s\in\Gamma_\pi}X_s-\sum_{s\in\Delta_\pi}X_s\label{eq:relC}
\end{equation}
as divisors on $\wt C$. Both sides of~(\ref{eq:relC}) are 
effective because $X_s$ and $Y_s$ are distinct for 
each $s\in\Gamma_\pi$. Moreover, since 
$X_s$ is not in the support of $W$ and $Y_s$ is not in the support of $W'$ 
for each $s\in\Gamma_\pi$, 
neither $X_s$ nor $Y_s$ is in the support of 
either side of~(\ref{eq:relC}). So the right-hand side of~(\ref{eq:olW}) 
is a relative Cartier divisor on $\wt C/S$. As both sides 
of~(\ref{eq:olW}) are relative Cartier divisors and agree on the generic 
fiber of $\wt\pi$, they agree everywhere.
\end{proof}

\section{Proof of the theorem}
\label{Sect6}

\begin{claim}\emph{Preliminaries.} 
Let $D\ol w\:{\cal O}_{\wt C}\to J^1_{\wt\pi}({\cal O}_{\wt C}(\ol W))$
be the derivative of the section $\ol w$ of ${\cal O}_{\wt C}(\ol W)$ 
given by $\ol W$. Since $\ol W$ is the disjoint union of subschemes 
isomorphic to $S$, the zero scheme $Z$ 
of $D\ol w$ is empty. Since $Z$ has the expected dimension, 
$[Z]=c_2(J^1_{\wt\pi}({\cal O}_{\wt C}(\ol W)))$. From exact 
sequence~(\ref{eq:short}) with $\wt\pi$ for $\pi$ and 
${\cal O}_{\wt C}(\ol W)$ for $\cal L$ we get
\[
c_2(J^1_{\wt\pi}({\cal O}_{\wt C}(\ol W)))=c_1({\cal O}_{\wt C}(\ol W))
c_1(\w_{\wt\pi}\otimes{\cal O}_{\wt C}(\ol W)).
\]

Let $K:=c_1(\w_{\wt\pi})$. Since $W$ is the zero scheme of a map
${\wt\pi}^*\det\wt\pi_*\cal L\to\w_{\wt\pi}\otimes{\cal L}^2$, we have 
$W=K+2c_1(\cal L)-\wt\pi^*c_1(\wt\pi_*\cal L)$. From~(\ref{eq:olW}) we get
\begin{equation}
\ol W=K+D-{\wt\pi}^*c_1(\wt\pi_*\cal L),\quad\text{where}\quad
D:=2c_1(\cal L)-\sum_{s\in\Xi_\pi}Y_s-2\sum_{s\in\Delta_\pi}Y_s.\label{eq:P1}
\end{equation}
By (\ref{eq:Lg-1}), there is an invertible sheaf $N$ on $S$ such that
\begin{equation}
(g-1)c_1(\cal L)=K+\sum_{s\in\Xi_\pi}b_sY_s+
\sum_{s\in\Delta_\pi}(2b_s-1)Y_s+\wt\pi^*c_1(N).\label{eq:P2}
\end{equation}
\end{claim}

\begin{claim}\emph{Computing $c_1(\wt\pi_*\cal L)$.} By definition, 
$\cal L={\cal O}_{\wt C}(2\Sigma)$. We 
claim that $(\wt\pi_*{\cal L}(s),{\cal L}_s)$ has no base points 
for any $s\in S$. This is well-known for $s\not\in\Gamma_\pi$. 
Let now $s\in\Gamma_\pi$. Since $\ol M_{0,2g+2}$ parameterizes 
\emph{stable} $2g+2$-pointed rational curves, 
there is an integer $i>1$ such that $\Sigma_{i,s}\subseteq X_s$. As 
${\cal O}_{\wt C}(2\Sigma_i)$ and $\cal L$ agree on the generic 
fiber of $\wt\pi$, and restrict to sheaves of equal degree on the 
irreducible components of $\wt C_s$, they agree on a neighborhood of 
$\wt C_s$. So $2\Sigma_s$ and $2\Sigma_{i,s}$ are divisors 
of the linear system $(\wt\pi_*{\cal L}(s),{\cal L}_s)$, which has 
therefore no base points.

It follows from our claim above that the natural sequence below is exact.
\[
0\to\wt\pi_*{\cal O}_{\wt C}(\Sigma)\to\wt\pi_*{\cal O}_{\wt C}(2\Sigma)\to
\wt\pi_*{\cal O}_{\wt C}(2\Sigma)|\Sigma\to 0.
\] 
So $\wt\pi_*{\cal O}_{\wt C}(\Sigma)$ is invertible on $S$. Now, 
the natural inclusion ${\cal O}_{\wt C}\to\cal O_{\wt C}(\Sigma)$ pushes 
down to a map 
${\cal O}_S\to\wt\pi_*{\cal O}_{\wt C}(\Sigma)$ whose restriction to 
each point of $S$ is injective. Since $\wt\pi_*{\cal O}_{\wt C}(\Sigma)$ is 
invertible, ${\cal O}_S=\wt\pi_*{\cal O}_{\wt C}(\Sigma)$. 
So, it follows from the above exact sequence that
\[
c_1(\wt\pi_*\cal L)=c_1(\wt\pi_*{\cal O}_{\wt C}(2\Sigma)|\Sigma)=
\wt\pi_*(2\Sigma^2).
\]
Since $c_1(\cal L)=2\Sigma$, it follows that
\begin{equation}
\wt\pi_*(c_1(\cal L)^2)=2c_1(\wt\pi_*\cal L).\label{eq:c1M}
\end{equation}
\end{claim}

\begin{claim}\emph{Conclusion.} For each $s\in\Gamma_\pi$, we have
\begin{equation}
\wt\pi_*(KY_s)=
\begin{cases}
2b_s[s]&\text{if $s\in\Xi_\pi$},\\
(2b_s-1)[s]&\text{if $s\in\Delta_\pi$}.
\end{cases}\label{eq:P5}
\end{equation}
Now, since $X_s+Y_s=\wt\pi^*[s]$, using the projection formula we 
get
\begin{equation}
\wt\pi_*Y_s^2=
\begin{cases}
-2[s]&\text{if $s\in\Xi_\pi$},\\
-[s]&\text{if $s\in\Delta_\pi$}.
\end{cases}\label{eq:P6}
\end{equation}
Let $\kappa_{1,\wt\pi}:=\wt\pi_*K^2$. 
Since $[Z]=\ol W(K+\ol W)$, it follows from~(\ref{eq:P1}) that
\begin{equation}
[Z]=(K+D-{\wt\pi}^*c_1(\wt\pi_*\cal L))(2K+D-{\wt\pi}^*c_1(\wt\pi_*\cal L)).
\label{eq:Z}
\end{equation}
Expanding, we get
\[
[Z]=2K^2+3KD+D^2-3K{\wt\pi}^*c_1(\wt\pi_*\cal L)
-2D{\wt\pi}^*c_1(\wt\pi_*\cal L)+({\wt\pi}^*c_1(\wt\pi_*\cal L))^2.
\]
Pushing down to $S$, and using the projection formula, we get
\[
\wt\pi_*[Z]=2\kappa_{1,\wt\pi}+\wt\pi_*(3KD+D^2)-6(g-1)c_1(\wt\pi_*\cal L)
-2\wt\pi_*Dc_1(\wt\pi_*\cal L).
\]
Now, $c_1(\cal L)Y_s=0$ for every $s\in\Gamma_\pi$, because 
$\Sigma_s\subseteq X_s$. Thus, using the expression for $D$ 
in~(\ref{eq:P1}), we get
\begin{align*}
\wt\pi_*[Z]=&2\kappa_{1,\wt\pi}+6\wt\pi_*(Kc_1(\cal L))
-6(g-1)c_1(\wt\pi_*\cal L)+4\wt\pi_*(c_1(\cal L)^2)-8c_1(\wt\pi_*\cal L)\\
&-\sum_{s\in\Xi_\pi}\wt\pi_*(3KY_s-Y_s^2)-
\sum_{s\in\Delta_\pi}\wt\pi_*(6KY_s-4Y_s^2).
\end{align*}
Using~(\ref{eq:c1M}),~(\ref{eq:P5})~and~(\ref{eq:P6}), 
the right-hand side of the above expression becomes
\[
2\kappa_{1,\wt\pi}+6\wt\pi_*(Kc_1(\cal L))-3(g-1)\wt\pi_*(c_1(\cal L)^2)
-\sum_{s\in\Xi_\pi}(6b_s+2)[s]-\sum_{s\in\Delta_\pi}(12b_s-2)[s].
\]
Using~(\ref{eq:P2}), we obtain
\[
\wt\pi_*[Z]=2\kappa_{1,\wt\pi}+3\wt\pi_*(Kc_1(\cal L))-6c_1(N)
-\sum_{s\in\Xi_\pi}(6b_s+2)[s]
-\sum_{s\in\Delta_\pi}(12b_s-2)[s].
\]
Multiplying by $g-1$, and using~(\ref{eq:P2}) again, 
the right-hand side of the above equation becomes
\[
(2g+1)\kappa_{1,\wt\pi}+
\sum_{x\in\Xi_\pi}(6b_s(b_s-g+1)-2(g-1))[s]
+\sum_{s\in\Delta_\pi}(12b_s(b_s-g)+2g+1)[s].
\]
Since $a_s+b_s=g-1$ if $s\in\Xi_\pi$ and $a_s+b_s=g$ if $s\in\Delta_\pi$, 
we get
\[
(g-1)\wt\pi_*[Z]=
(2g+1)\kappa_{1,\wt\pi}-\sum_{s\in\Xi_\pi}(6a_sb_s+2(g-1))[s]
-\sum_{s\in\Delta_\pi}(12a_sb_s-2g-1)[s].
\]
Now, by Grothendieck-Riemann-Roch,
\[
\kappa_{1,\wt\pi}=12\lambda_{\wt\pi}-2\sum_{s\in\Xi_\pi}[s]
-\sum_{s\in\Delta_\pi}[s].
\]
Moreover, since $\wt\pi_*\w_{\wt\pi}=\pi_*\w_\pi$, we have 
$\lambda_{\wt\pi}=\lambda_\pi$. Thus
\begin{align*}
(g-1)\wt\pi_*[Z]=&12(2g+1)\lambda_\pi-6\sum_{x\in\Xi_\pi}(a_sb_s+g)[s]-
12\sum_{s\in\Delta_\pi}a_sb_s[s]\\
=&12(2g+1)\lambda_\pi-6\sum_{i=0}^{[(g-1)/2]}(i+1)(g-i)[\Xi_{i,\pi}]-
12\sum_{j=1}^{[g/2]}j(g-j)[\Delta_{j,\pi}].
\end{align*}
Using (\ref{eq:DeltaXipi1})~and~(\ref{eq:DeltaXipi2}) we get
\[
(g-1)\wt\pi_*[Z]=12(2g+1)\lambda_{\pi}-3g\xi_{0,\pi}
-6\sum_{i=1}^{[(g-1)/2]}(i+1)(g-i)\xi_{\pi,i}-
12\sum_{j=1}^{[g/2]}j(g-j)\delta_{\pi,j}.
\]
Then $(g-1)\wt\pi_*[Z]$ is a multiple of the left-hand side 
of~(\ref{eq:CHpi}). Since $\wt\pi_*[Z]=0$, we get~(\ref{eq:CHpi}), thus 
proving Cornalba's and Harris' relation~(\ref{eq:CH}).
\end{claim}

\section{A variant}
\label{Sect7}

There is yet another way to obtain relation~(\ref{eq:CH}), 
still inspired by Faber's general philosophy of writing down the class 
of the empty set. We present here just a rough account of this different 
approach, as we intend to give details in a forthcoming article.

Keep the notation used so far in the article. 
The particular empty set we will consider is the closure $\ol R_{1,2}$ 
of the locus $R_{1,2}\subseteq\wt C\times_S\wt C$ of 
pairs $(P,Q)$ of points of a smooth fiber $\wt C_s$ of 
$\wt\pi$ for which $\cal O_{\wt C_s}(P+2Q)\cong\cal L|\wt C_s$. This 
set is clearly empty, as $\cal L$ has relative degree 2 over $S$. 
On the other hand, we can express the class of 
$\ol R_{1,2}$ in $A^2(\wt C\times_S\wt C)$ 
as the product,
\begin{equation}
(L_1+L_2-\Delta-\rho^*c_1(\wt\pi_*\cal L)-F)
(K_2+L_1+L_2-\Delta-\rho^*c_1(\wt\pi_*\cal L)-F),\label{eq:variant}
\end{equation}
where $\rho\:\wt C\times_S\wt C\to S$ is the natural map, 
$K_i:=p_i^*c_1(\w_\pi)$ and $L_i:=p_i^*c_1(\cal L)$ for $i=1,2$, 
where $p_1$ and $p_2$ are the projections of $\wt C\times_S\wt C$ onto the 
indicated factors, $\Delta$ is the diagonal of $\wt C\times_S\wt C$ and 
$F$ is a correction term supported in the singular fibers of $\rho$. 
The computation of $F$ is similar to that which leads to 
relation~(\ref{eq:olW}). 

Being~(\ref{eq:variant}) zero, we obtain a relation 
in $A^2(\wt C\times_S\wt C)$. By reason of dimension, this relation does not 
yield immediately a relation in $\text{Pic}(S)$. However, using 
Faber's recipe in~\cite{Faber3}, ``once zero always zero'', 
we multiply~(\ref{eq:variant}) by divisor classes on 
$\wt C\times_S\wt C$, and then push down the ensuing products to $S$ in 
order to get relations in $\text{Pic}(S)$. There are two natural 
divisor classes on $\wt C\times_S\wt C$: the diagonal 
$\Delta$ and $K_1$. Hence we get two relations in $\text{Pic}(S)$:
\begin{align}
\rho_*((L_1+L_2-\Delta-\rho^*c_1(\wt\pi_*\cal L)-F)
(K_2+L_1+L_2-\Delta-\rho^*c_1(\wt\pi_*\cal L)-F)&\Delta)=0,\\ 
\rho_*((L_1+L_2-\Delta-\rho^*c_1(\wt\pi_*\cal L)-F)
(K_2+L_1+L_2-\Delta-\rho^*c_1(\wt\pi_*\cal L)-F)&K_1)=0.
\end{align}
The left-hand side of~(25) is simply $\wt\pi_*[Z]$, for $[Z]$ expressed as 
in~(\ref{eq:Z}). We use~(26) to find an expression for 
$c_1(\wt\pi_*\cal L)$ which, when 
used in~(25), yields~(\ref{eq:CHpi}) after a few computations, as in 
Subsection~6.3.

The main difficulty in the above alternative proof of~(\ref{eq:CH}) 
is that $\wt C\times_S\wt C$ is not smooth. More precisely, $\Delta$ and 
the irreducible components of the singular fibers of $\rho$ are singular in 
$\wt C\times_S\wt C$, and hence equations~(24)-(26) don't make immediate 
sense. This difficulty can be overcome because the singular 
locus of $\wt C\times_S\wt C$ is finite, of codimension 3.

\vspace{5mm}
\begin{flushright}\begin{minipage}[r]{2.5in}
\baselineskip=12pt
\midrm Instituto de Matem\'atica Pura e Aplicada\\
Estrada Dona Castorina, 110\\
22460-320 Rio de Janeiro RJ (BRAZIL)\\
esteves@impa.br
\end{minipage}
\end{flushright}

\vspace{5mm}
\begin{flushright}\begin{minipage}[r]{2.5in}
\baselineskip=12pt
\midrm Departamento de Matem\'atica\\
Universidade Federal de Pernambuco\\
Av. Prof. Luiz Freire S/N\\
Cidade Universit\'aria\\
50540-740 Pernambuco PE (BRAZIL)\\
taggo@dmat.ufpe.br
\end{minipage}
\end{flushright}

\vspace{5mm}
\begin{flushright}\begin{minipage}[r]{2.5in}
\baselineskip=12pt
\midrm Dipartimento di Matematica\\
Politecnico di Torino\\
C.so Duca degli Abruzzi 24\\
10129 Torino (ITALY)\\
lgatto@polito.it
\end{minipage}
\end{flushright}

\end{document}